УДК 511.72+517.5+517.13

**Symon Serbenyuk**
*e-mail:simon6@ukr.net*

# GENERALIZED SHIFT OPERATOR OF CERTAIN ENCODINGS OF REAL NUMBERS

*The present article is devoted to the investigation of some properties of the generalized shift operator of numbers represented in terms of numeral systems with a variable alphabet.*

1. **Introduction.** Let $N_B$ be a fixed subset of positive integers, $B=(b_n)$ be a fixed increasing sequence of all elements of $N_B$, $\rho_0=0$,

$$\rho_n=\begin{cases}1 \ whenever \ n\in N_B\\ 2 \ whenever \ n\notin N_B\end{cases},$$

and

$$\tilde{Q}'_{N_B}=\begin{pmatrix} (-1)^{\rho_1}q_{0,1} & (-1)^{\rho_1+\rho_2}q_{0,2} & \cdots & (-1)^{\rho_{n-1}+\rho_n}q_{0,n} & \cdots\\ (-1)^{\rho_1}q_{1,1} & (-1)^{\rho_1+\rho_2}q_{1,2} & \cdots & (-1)^{\rho_{n-1}+\rho_n}q_{1,n} & \cdots\\ \vdots & \vdots & \ddots & \vdots & \cdots\\ (-1)^{\rho_1}q_{m_1-1,1} & (-1)^{\rho_1+\rho_2}q_{m_2-2,2} & \cdots & (-1)^{\rho_{n-1}+\rho_n}q_{m_n,n} & \cdots\\ & (-1)^{\rho_1+\rho_2}q_{m_2-1,2} & \cdots & & \cdots\\ & (-1)^{\rho_1+\rho_2}q_{m_2,2} & \cdots & & \cdots \end{pmatrix}$$

be a fixed matrix. Here $i=\overline{0,m_n}$, $m_n\in N\cup\{0,+\infty\}$, $n=1,2,...,$ and elements $q_{i,n}$ satisfy the following system of conditions:

$$\begin{cases}1.\ 0<q_{i,j}\in R\\ 2.\ \forall j\in N:\sum_{i=0}^{m_j}q_{i,j}=1\\ 3.\ \forall(i_j),\ i_j\in N\cup\{0\}:\prod_{j=1}^{\infty}q_{i_j,j}=0.\end{cases}$$

So, any number from a certain interval $(a',a'')$ can be represented by the following way (see [14, 9]):

$$x=\Delta^{\tilde{Q}'_{N_B}}_{i_1i_2...i_n...}=(-1)^{\rho_1}a_{i_1(x),1}+\sum_{k=2}^{\infty}\left[(-1)^{\rho_k}a_{i_k(x),k}\prod_{j=1}^{k-1}q_{i_j(x),j}\right], \quad (1)$$

where

$$a_{i_k(x),k}=\begin{cases}\sum_{i=0}^{i_k-1}q_{i,k} & whenever\ i_k\neq 0\\ 0 & whenever\ i_k=0.\end{cases}$$

Note that representation (1) is:

- the nega-$\tilde{Q}$-representation ([10]) whenever a set $N_B$ is the set of all odd positive integers, i.e.,

$$x=\Delta^{-\tilde{Q}}_{i_1i_2...i_n...}=-a_{i_1,1}+\sum_{n=2}^{\infty}\left[(-1)^n a_{i_n,n}\prod_{j=1}^{n-1}q_{i_j,j}\right];$$

- the representation by a positive Cantor series (such series were introduced in [1]) whenever $N_B$ is empty and the condition $q_{i,n} = \frac{1}{q_n}$ holds for all $n \in N$, where $Q = (q_n)$ is a fixed sequence of positive integers and $q_n > 1$:

$$x = \Delta^Q_{i_1 i_2 \ldots i_n \ldots} = \sum_{n=1}^{\infty} \frac{i_n}{q_1 q_2 \ldots q_n}, \quad i_n \in \{0,1,\ldots,q_n - 1\}; \tag{2}$$

- the representation by an alternating Cantor series ([11, 5,6]) whenever the last-mentioned conditions hold but $N_B$ is the set of all odd positive integers, i.e.,

$$x = \Delta^Q_{i_1 i_2 \ldots i_n \ldots} = \sum_{n=1}^{\infty} \frac{(-1)^n i_n}{q_1 q_2 \ldots q_n}, \quad i_n \in \{0,1,\ldots,q_n - 1\}; \tag{3}$$

- the q-ary representation ([3]) of real numbers whenever $q_{i,n} = \frac{1}{q}$ for all $n \in N$, where $q > 1$ is a fixed positive integer, and $N_B$ is empty:

$$x = \Delta^q_{i_1 i_2 \ldots i_n \ldots} = \sum_{n=1}^{\infty} \frac{i_n}{q^n}, \quad i_n \in \{0,1,\ldots,q - 1\};$$

- is the nega-q-ary representation ([2]) whenever $q_{i,n} = \frac{1}{q}$ for all $n \in N$, where $q > 1$ is a fixed positive integer, and $N_B$ is the set of all odd positive integers:

$$x = \Delta^{-q}_{i_1 i_2 \ldots i_n \ldots} = \sum_{n=1}^{\infty} \frac{i_n}{(-q)^n}, \quad i_n \in \{0,1,\ldots,q - 1\}.$$

In 2013, a notion of the generalized shift operator was defined for alternating Cantor series (see the working paper [6], the presentation [5], and the corresponding published paper [11]). The generalized shift operator is following in terms of alternating Cantor series:

$$\sigma_m\left(\sum_{n=1}^{\infty} \frac{(-1)^m i_n}{q_1 q_2 \ldots q_n}\right) = -\frac{i_1}{q_1} + \frac{i_2}{q_1 q_2} - \frac{i_3}{q_1 q_2 q_3} + \ldots + \frac{(-1)^{m-1} i_{m-1}}{q_1 q_2 \ldots q_{m-1}} + \frac{(-1)^m i_{m+1}}{q_1 q_2 \ldots q_{m-1} q_{m+1}} + \frac{(-1)^{m+1} i_{m+2}}{q_1 q_2 \ldots q_{m-1} q_{m+1} q_{m+2}} + \ldots$$

The idea includes the following: any number from a certain interval can be represented by two fixed sequences $(q_n)$ and $(i_n)$. The generalized shift operator maps the preimage into a number represented by the following two sequences $(q_1, q_2, \ldots, q_{m-1}, q_{m+1}, \ldots)$ and $(i_1, i_2, \ldots, i_{m-1}, i_{m+1}, \ldots)$. In terms of certain encodings of real numbers, this number can belong to another interval.

In the present article, the generalized shift operator is investigated for different expansions of real numbers (the main attention is given to numeral systems with variable alphabets).

Let us remark that some numeral system is a numeral system with a variable alphabet whenever there exist at least two numbers $k$ and $l$ such that the condition $A_k \neq A_l$ holds for the representation $\Delta_{i_1 i_2 \ldots i_n \ldots}$ of numbers in terms of this numeral system, where $i_k \in A_k$ and $i_l \in A_l$, as well as $k \neq l$.

2. **Cantor series.** Now let us consider numbers represented in terms of different Cantor series. We begin with positive Cantor series.

Suppose a number $x \in [0,1]$ represented by positive Cantor series (1). Then

$$\sigma_m(x) = \sum_{k=1}^{m-1} \frac{i_k}{q_1 q_2 \ldots q_k} + \sum_{l=m+1}^{\infty} \frac{i_l}{q_1 q_2 \ldots q_{m-1} q_{m+1} \ldots q_l}.$$

Denote by $\varsigma_{m+1}$ the sum $\sum_{l=m+1}^{\infty} \frac{i_l}{q q \ldots q_{m-1} q_{m+1} \ldots q_l}$ and by $\vartheta_{m-1}$ the sum $\sum_{k=1}^{m-1} \frac{i_k}{q_1 q_2 \ldots q_k}$. Then $\varsigma_{m+1} = q_m(x - \vartheta_m)$ and

$$\sigma_m(x) = q_m x - (q_m - 1)\vartheta_{m-1} - \frac{i_m}{q_1 q_2 \ldots q_{m-1}}. \tag{4}$$

Let us remark that

$$\sigma(x)=\sigma_1(x)=\sum_{n=2}^{\infty}\frac{i_n}{q_1q_2\ldots q_n}=q_1\Delta^Q_{0i_2i_3\ldots i_n}$$

and, where $\sigma$ is the shit operator,

$$\sigma^m(x)=\sum_{k=m+1}^{\infty}\frac{i_k}{q_{m+1}q_{m+2}\ldots q_k}=q_1q_2\ldots q_m\Delta^Q_{\underbrace{0\ldots0}_{m}i_{m+1}i_{m+2}\ldots}.$$

Whence

$$\sigma_m(x)=\Delta^Q_{i_1i_2\ldots i_{m-1}000\ldots}+q_m\Delta^Q_{\underbrace{0\ldots0}_{m}i_{m+1}i_{m+2}\ldots}=\Delta^Q_{i_1i_2\ldots i_{m-1}0i_{m+1}i_{m+2}\ldots}+(q_m-1)\Delta^Q_{\underbrace{0\ldots0}_{m}i_{m+1}i_{m+2}\ldots}.$$

Note that for positive Cantor series numbers of the form

$$\Delta^Q_{i_1i_2\ldots i_{n-1}i_n000\ldots}=\Delta^Q_{i_1i_2\ldots i_{n-1}[i_n-1][q_{n+1}-1][q_{n+2}-1][q_{n+3}-1]\ldots}$$

are called *Q-rational*. The other numbers in $[0,1]$ are called *Q-irrational*.

**Theorem.** *The generalized shift operator has the following properties:*

- $\sigma\circ\sigma_2^m=\sigma^{m+1}(x).$
- *Suppose* $(k_n)$ *is a sequence of positive integers such that* $k_n=k_{n-1}+1,\ n=2,3,\ldots$. *Then*
$$\sigma^{k_1+1}\circ\sigma_{k_n}\circ\sigma_{k_{n-1}}\circ\ldots\circ\sigma_{k_1}(x)=\sigma^{k_n}(x).$$
- *Suppose* $(k_n)$ *is an arbitrary finite subsequence of positive integers. Then*
$$\sigma^{k_n-n}\circ\sigma_{k_n}\circ\sigma_{k_{n-1}}\circ\ldots\circ\sigma_{k_1}(x)=\sigma^{k_n}(x).$$
- *The mapping* $\sigma_m$ *is continuous at any Q-irrational point and is continuous at any Q-rational point*
$$x_0=x_0^{(1)}=\Delta^Q_{i_1i_2\ldots i_{n-1}i_n000\ldots}=\Delta^Q_{i_1i_2\ldots i_{n-1}[i_n-1][q_{n+1}-1][q_{n+2}-1][q_{n+3}-1]\ldots}=x_0^{(2)}$$
*whenever* $m\neq n$. *If* $m=n,$ *then* $x_0$ *is a point of discontinuity of* $\sigma_m$.
- *If the mapping* $\sigma_m$ *has a derivative at the point* $x=\Delta^Q_{i_1i_2\ldots i_n\ldots}$, *then* $(\sigma_m)'=q_m$.
- $$x-\sigma_m(x)=\frac{i_m}{q_1q_2\ldots q_m}+\frac{\sigma^m(x)}{q_1q_2\ldots q_m}(1-q_m).$$

*Proof.* All properties follow from the definition of $\sigma_m$ and equality (4).

Let us consider the set $\Delta^Q_{c_1c_2\ldots c_n}$ of all numbers in whose representations by Cantor series the first $n$ digits, i.e., $c_1,c_2,\ldots,c_n$, are fixed but other digits are arbitrary from $A_k=\{0,1,\ldots,q_k-1\},\ k=n+1,n+2,\ldots$. It is easy to see that $\lim\limits_{x\to x_0}\sigma_m(x)=x_0$ for all numbers from $\Delta^Q_{c_1c_2\ldots c_n}$ whenever $m\neq n$. If $m=n$, then

$$\lim_{x\to x_0+0}\sigma_m(x)=\sigma_m(x_0^1),\ \lim_{x\to x_0-0}\sigma_m(x)=\sigma_m(x_0^2),\text{ and }\sigma_m(x_0^1)-\sigma_m(x_0^2)=-\frac{1}{q_1q_2\ldots q_{m-1}}.$$

Also,

$$x-\sigma_m(x)=\vartheta_m+\frac{\sigma^m(x)}{q_1q_2\ldots q_m}-\vartheta_{m-1}-\varsigma_{m+1}=\frac{i_m}{q_1\ldots q_m}+\frac{\sigma^m(x)}{q_1q_2\ldots q_m}-\frac{\sigma^m(x)}{q_1q_2\ldots q_{m-1}}=\frac{i_m}{q_1\ldots q_m}+\frac{\sigma^m(x)}{q_1q_2\ldots q_m}(1-q_m).$$

□

Let us consider expansion (3). In this case,

$$\sigma_m(x)=\sigma_m\left(\sum_{n=1}^{\infty}\frac{(-1)^n i_n}{q_1q_2\ldots q_n}\right)=\sum_{k=1}^{m-1}\frac{(-1)^k i_k}{q_1\ldots q_k}+\sum_{l=m+1}^{\infty}\frac{(-1)^{l-1}i_l}{q_1\ldots q_{m-1}q_{m+1}\ldots q_l}$$

$$=-q_mx+(1+q_m)\sum_{k=1}^{m-1}\frac{(-1)^k i_k}{q_1q_2\ldots q_k}+\frac{(-1)^m i_m}{q_1\ldots q_{m-1}}.$$

So, $\sigma_m$ is a piecewise linear function since $\sum\limits_{k=1}^{m}\frac{(-1)^k c_k}{q_1\ldots q_k}$ is constant for the set $\Delta^{-Q}_{c_1c_2\ldots c_m}$.

In the general case, we have

$$x=\Delta^{(\pm Q,N_B)}_{i_1i_2\dots i_n\dots}=\sum_{n=1}^{\infty}\frac{(-1)^{\rho_n}i_n}{q_1q_2\dots q_n},\ i_n\in A_n,$$

and

$$\sigma_m(x)=\frac{(-1)^{\rho_1}i_1}{q_1}+\frac{(-1)^{\rho_2}i_2}{q_1q_2}+\dots+\frac{(-1)^{\rho_{m-1}}i_{m-1}}{q_1\dots q_{m-1}}+\frac{(-1)^{\rho_{m+1}}i_{m+1}}{q_1\dots q_{m-1}q_{m+1}}+\dots=\sum_{k=1}^{m-1}\frac{(-1)^{\rho_k}i_k}{q_1\dots q_k}+\sum_{l=m+1}^{\infty}\frac{(-1)^{\rho_l}i_l}{q_1\dots q_{m-1}q_{m+1}\dots q_l}.$$

Since

$$x=\sum_{k=1}^{m-1}\frac{(-1)^{\rho_k}i_k}{q_1\dots q_k}+\frac{(-1)^{\rho_m}i_m}{q_1q_2\dots q_m}+\frac{1}{q_m}\sum_{l=m+1}^{\infty}\frac{(-1)^{\rho_l}i_l}{q_1\dots q_{m-1}q_{m+1}\dots q_l},$$

we get

$$\sigma_m(x)=q_mx+(1-q_m)\sum_{k=1}^{m-1}\frac{(-1)^{\rho_k}i_k}{q_1\dots q_k}-\frac{(-1)^{\rho_m}i_m}{q_1q_2\dots q_{m-1}}.$$

By analogy, only for the case when $m=n$, points of discontinuity of $\sigma_m$ are only quasi-Q-rational points which are points of the following form (for example, quasi-q-ary rational numbers were described in [9]):

- if $n,n+1\in N_B$, then $\Delta^{(\pm Q,N_B)}_{i_1i_2\dots i_{n-1}[q_n-1-i_n][q_{n+1}-1]\beta_{n+2}\beta_{n+3}\dots}=\Delta^{(\pm Q,N_B)}_{i_1i_2\dots i_{n-1}[q_n-i_n]0\gamma_{n+2}\gamma_{n+3}\dots}$;
- if $n\in N_B,(n+1)\notin N_B$, then $\Delta^{(\pm Q,N_B)}_{i_1i_2\dots i_{n-1}[q_n-1-i_n]0\beta_{n+2}\beta_{n+3}\dots}=\Delta^{(\pm Q,N_B)}_{i_1i_2\dots i_{n-1}[q_n-i_n][q_{n+1}-1]\gamma_{n+2}\gamma_{n+3}\dots}$;
- if $n\notin N_B,(n+1)\in N_B$, then $\Delta^{(\pm Q,N_B)}_{i_1i_2\dots i_{n-1}i_n[q_{n+1}-1]\beta_{n+2}\beta_{n+3}\dots}=\Delta^{(\pm Q,N_B)}_{i_1i_2\dots i_{n-1}[i_n-1]0\gamma_{n+2}\gamma_{n+3}\dots}$;
- if $n\notin N_B,(n+1)\notin N_B$, then $\Delta^{(\pm Q,N_B)}_{i_1i_2\dots i_{n-1}i_n0\beta_{n+2}\beta_{n+3}\dots}=\Delta^{(\pm Q,N_B)}_{i_1i_2\dots i_{n-1}[i_n-1][q_{n+1}-1]\gamma_{n+2}\gamma_{n+3}\dots}$.

Here $\beta_n=\begin{cases}0 & whenever\ n\notin N_B\\ q_n-1 & whenever\ n\in N_B\end{cases}$ and $\gamma_n=\begin{cases}q_n-1 & whenever\ n\notin N_B\\ 0 & whenever\ n\in N_B.\end{cases}$

In other cases, $\sigma_m$ is continuous and $(\sigma_m)'=q_m$.

**3. $\tilde{Q}-$ expansions.** Let us consider a positive $\tilde{Q}-$ expansion. Suppose that

$$x=\Delta^{\tilde{Q}}_{i_1i_2\dots i_n\dots}=a_{i_1,1}+\sum_{n=2}^{\infty}\left(a_{i_n,n}\prod_{j=1}^{n-1}q_{i_j,j}\right).$$

Then

$$\sigma_m(x)=a_{i_1,1}+\sum_{n=2}^{m-1}\left(a_{i_n,n}\prod_{j=1}^{n-1}q_{i_j,j}\right)+\sum_{l=m+1}^{\infty}\left(a_{i_l,l}\prod_{\substack{j=1,\\ j\neq m}}^{l-1}q_{i_j,j}\right).$$

Since

$$x=a_{i_1,1}+\sum_{n=2}^{m-1}\left(a_{i_n,n}\prod_{j=1}^{n-1}q_{i_j,j}\right)+a_{i_m,m}\prod_{j=1}^{m-1}q_{i_j,j}+q_{i_m,m}\sum_{l=m+1}^{\infty}\left(a_{i_l,l}\prod_{\substack{j=1,\\ j\neq m}}^{l-1}q_{i_j,j}\right),$$

we have

$$\sigma_m(x)=a_{i_1,1}+\sum_{n=2}^{m-1}\left(a_{i_n,n}\prod_{j=1}^{n-1}q_{i_j,j}\right)+\frac{x}{q_{i_m,m}}-\frac{1}{q_{i_m,m}}\left(a_{i_1,1}+\sum_{n=2}^{m-1}\left(a_{i_n,n}\prod_{j=1}^{n-1}q_{i_j,j}\right)+a_{i_m,m}\prod_{j=1}^{m-1}q_{i_j,j}\right)$$

$$=\frac{x}{q_{i_m,m}}-\frac{a_{i_m,m}\prod_{j=1}^{m-1}q_{i_j,j}}{q_{i_m,m}}+\left(1-\frac{1}{q_{i_m,m}}\right)\left(a_{i_1,1}+\sum_{n=2}^{m-1}\left(a_{i_n,n}\prod_{j=1}^{n-1}q_{i_j,j}\right)\right).$$

In addition, note that

$$\sigma^m(x) = \sigma_1^m(x) = \frac{1}{q_{0,1}q_{0,2}\ldots q_{0,m}}\Delta^{\tilde{Q}}_{\underbrace{0\ldots0}_{m}i_{m+1}i_{m+2}\ldots},$$

$$x = \Delta^{\tilde{Q}}_{i_1i_2\ldots i_n000\ldots} + \sigma^m(x)\prod_{j=1}^{m} q_{i_j,j},$$

and

$$\sigma_m(x) = \Delta^{\tilde{Q}}_{i_1i_2\ldots i_{m-1}000\ldots} + \frac{\prod_{j=1}^{m-1} q_{i_j,j}}{\prod_{t=1}^{m} q_{0,t}}\Delta^{\tilde{Q}}_{\underbrace{0\ldots0}_{m}i_{m+1}i_{m+2}\ldots}.$$

For the general case

$$x = \Delta^{(\pm\tilde{Q},N_B)}_{i_1i_2\ldots i_n\ldots} = (-1)^{\rho_1}a_{i_1,1} + \sum_{n=2}^{\infty}\left((-1)^{\rho_n}a_{i_n,n}\prod_{j=1}^{n-1}q_{i_j,j}\right),$$

we have

$$\sigma_m(x) = (-1)^{\rho_1}a_{i_1,1} + \sum_{n=2}^{m-1}\left((-1)^{\rho_n}a_{i_n,n}\prod_{j=1}^{n-1}q_{i_j,j}\right) + \sum_{l=m+1}^{\infty}\left((-1)^{\rho_l}a_{i_l,l}\prod_{\substack{j=1,\\ j\ne m}}^{l-1}q_{i_j,j}\right).$$

and

$$\sigma_m(x) = \frac{x}{q_{i_m,m}} - \frac{(-1)^{\rho_m}a_{i_m,m}\prod_{j=1}^{m-1}q_{i_j,j}}{q_{i_m,m}} + \left(1 - \frac{1}{q_{i_m,m}}\right)\left((-1)^{\rho_1}a_{i_1,1} + \sum_{n=2}^{m-1}\left((-1)^{\rho_n}a_{i_n,n}\prod_{j=1}^{n-1}q_{i_j,j}\right)\right).$$

From considered relationships follows the following statement.

**Theorem.** *For any case of a sign variable $\tilde{Q}$ - expansion, the mapping $\sigma_m$ is a piecewise linear function.*

**4. Conclusions.** In the present article, the attention was given to the generalized shift operator in terms of certain numeral systems with variable alphabets. Some relationships between $x, \sigma^m(x)$, and $\sigma_m(x)$ are considered in terms of some encodings of real numbers. The general property of $\sigma_m$ is a fact that this mapping is a piecewise linear function. Note that for numeral systems with variable alphabets, the following is true: $\sigma_m(\Delta_{i_1i_2\ldots i_n\ldots}) \ne \Delta_{i_1i_2\ldots i_{m-1}i_{m+1}\ldots}$. That is, in our case, this mapping maps $x = \Delta_{i_1i_2\ldots i_n\ldots}$ into a number represented in terms of the other numeral system. However, for some numeral systems with a constant alphabet (i.e., the q-ary or nega-q-ary representation) the condition $\sigma_m(\Delta_{i_1i_2\ldots i_n\ldots}) = \Delta_{i_1i_2\ldots i_{m-1}i_{m+1}\ldots}$ holds. In addition, let us remark that the notion of the generalized shift operator is useful for further investigations. Some such investigations were described in the preprint [15]. For example, one can model generalizations of the singular Salem function (the singular Salem function was introduced in [4]) by the generalized shift operator (partial cases of such generalizations are generalized Salem functions investigated in [7, 8, 12, 13]). Next papers of the author of this article will be devoted to such investigations.

## УЗАГАЛЬНЕНИЙ ОПЕРАТОР ЗСУВУ ЦИФР (СИМВОЛІВ) ПЕВНИХ КОДУВАНЬ ДІЙСНИХ ЧИСЕЛ

*Стаття присвячена вивченню деяких властивостей узагальненого оператора зсуву цифр (символів) подання дійсних чисел в термінах систем числення зі змінним алфавітом.*

**С. Сербенюк**
*e-mail:simon6@ukr.net*

## ОБОБЩЕННЫЙ ОПЕРАТОР СДВИГА ЦИФР (СИМВОЛОВ) ОПРЕДЕЛЕННЫХ КОДИРОВОК ДЕЙСТВИТЕЛЬНЫХ ЧИСЕЛ

*Статья посвящена изучению свойств обобщенного оператора сдвига цифр (символов) представления действительных чисел в терминах систем счисления с переменным алфавитом.*

**С. Сербенюк**
*e-mail:simon6@ukr.net*